\documentclass[12pt]{amsart}
\usepackage[T1]{fontenc}
\usepackage{lmodern,amsmath,amssymb,mathtools,microtype}
\usepackage[a4paper,margin=23mm]{geometry}
\usepackage{xurl}
\usepackage[hidelinks]{hyperref}
\newtheorem{theorem}{Theorem}
\newtheorem{lemma}[theorem]{Lemma}
\newtheorem{proposition}[theorem]{Proposition}
\newtheorem{corollary}[theorem]{Corollary}
\newtheorem{externalinput}[theorem]{External input}
\theoremstyle{definition}

\theoremstyle{remark}
\newtheorem{remark}[theorem]{Remark}
\numberwithin{equation}{section}
\newcommand{\R}{\mathbb R}
\newcommand{\Z}{\mathbb Z}

\newcommand{\C}{\mathbb C}

\newcommand{\1}{\mathbf 1}
\newcommand{\eps}{\varepsilon}
\newcommand{\dd}{\,d}
\newcommand{\supp}{\operatorname{supp}}

\newcommand{\pn}{p_n}
\newcommand{\qn}{q_n}

\title[The square functions for cones and elliptic quadrics]{The square function estimates for cones in arbitrary dimensions and elliptic quadrics}
\author{Faruk Temur}
\date{November 3, 2026}
\hypersetup{pdftitle={Square function estimates for cones in arbitrary dimensions and elliptic quadrics},pdfauthor={Faruk Temur},pdfsubject={A deduction from the October 2026 local-smoothing and auxiliary-energy statements}}
\begin{document}
\begin{abstract}
We demonstrate that  the  square function estimate for the light cone in $\R^d$, $d\ge3$, at the exponent $p=2(d-1)/(d-2)$, with loss $[\log(2R)]^{5/p}$ can be deduced from  local-smoothing and auxiliary-energy statements in the October 25, 2026 OpenAI manuscript.  Exact affine sections preserve these logarithmic bounds and give corresponding estimates for spheres, ellipsoids, elliptic paraboloids, and compact patches of two-sheet elliptic hyperboloids in $\R^{d-1}$. 
\end{abstract}
\maketitle

\section{Introduction}
Let $n\ge2$, and consider the truncated cone
\[
 \Gamma_n=\{(\xi,\tau)\in\R^n\times\R:\tau=|\xi|,\quad1\le|\xi|\le2\}.
\]
Thus the ambient cone dimension is $d=n+1$. For $R\ge2$, put $h=R^{-1}$ and choose a maximal $\sqrt h$-separated set $\{\omega_\theta\}\subset S^{n-1}$. For a fixed $C_0\ge2$, define
\begin{equation}\label{eq:sectorintro}
 \theta=\left\{(\xi,\tau):1\le|\xi|\le2,\quad
 \left|\frac\xi{|\xi|}-\omega_\theta\right|\le C_0\sqrt h,
 \quad|\tau-|\xi||\le C_0h\right\}.
\end{equation}
These sectors have bounded overlap and lie in rectangular frequency planks of dimensions $1,(\sqrt h)^{n-1},h$. We give the smooth partition and the dual physical scales in Section~\ref{sec:target}.

\begin{theorem}[Cone square function estimate]\label{thm:main}
Assume the October local-smoothing and auxiliary-energy inputs stated in Section~\ref{sec:inputs}. Let
\[
 p_n=\frac{2n}{n-1}=\frac{2(d-1)}{d-2}.
\]
For every $R\ge2$ and every Schwartz family with $\supp\widehat f_\theta\subset\theta$,
\begin{equation}\label{eq:maintheorem}
 \left\|\sum_\theta f_\theta\right\|_{L^{p_n}(\R^{n+1})}
 \le C_{n,C_0}[\log(2R)]^{5/p_n}
 \left\|\left(\sum_\theta|f_\theta|^2\right)^{1/2}\right\|_{L^{p_n}(\R^{n+1})}.
\end{equation}
The constant is independent of the angular net and the family. Fixed positive radial bands and fixed sector enlargements are allowed, with a corresponding change of constant. Moreover, for every $2\le p<p_n$, put
\[
 \beta_n(p)=n\left(1-\frac2p\right)\left(\frac5{p_n}+1\right).
\]
Then
\begin{equation}\label{eq:target}
 \left\|\sum_\theta f_\theta\right\|_p
 \le C_{n,p,C_0}[\log(2R)]^{\beta_n(p)}
 \left\|\left(\sum_\theta|f_\theta|^2\right)^{1/2}\right\|_p.
\end{equation}
At $p=2$, $\beta_n(2)=0$ and the constant is uniform in $R$. For $2<p<p_n$, the assertion also uses the spatial plate maximal theorem stated in Section~\ref{sec:inputs}. In particular, for every $2\le p\le p_n$ and $\varepsilon>0$, these bounds imply the weaker estimate with $C_{n,p,\varepsilon,C_0}R^\varepsilon$ loss.
\end{theorem}

The pieces are arbitrary within their sectors. They need not be free waves and need not come from one prescribed multiplier partition. The right side is the norm of the pointwise square function.
The dimension in the elliptic consequence is one lower. For a smooth hypersurface $\Sigma\subset\R^n$, a canonical $\delta$-cap means a piece of its $O(\delta)$ neighborhood with tangential diameter $O(\sqrt\delta)$. The centers are separated at scale $\sqrt\delta$ and the fixed enlargements have bounded overlap. Section~\ref{sec:elliptic} gives explicit definitions for each surface below.

\begin{theorem}[Elliptic quadric square function estimate]\label{thm:ellipticintro}\label{cor:elliptic}
Under the hypotheses of Theorem~\ref{thm:main}, let $\Sigma\subset\R^n$ be a fixed sphere or ellipsoid, a fixed compact patch of an elliptic paraboloid, or a fixed compact patch of one sheet of a two-sheet elliptic hyperboloid. For $0<\delta\le1/2$, canonical $\delta$-caps $\vartheta$, and arbitrary Schwartz functions with $\supp\widehat g_\vartheta\subset\vartheta$,
\begin{equation}\label{eq:ellipticintro}
 \left\|\sum_\vartheta g_\vartheta\right\|_{L^{p_n}(\R^n)}
 \le C_\Sigma[\log(2/\delta)]^{5/p_n}
 \left\|\left(\sum_\vartheta|g_\vartheta|^2\right)^{1/2}\right\|_{L^{p_n}(\R^n)},
 \qquad p_n=\frac{2n}{n-1}.
\end{equation}
For every $2\le p<p_n$, the same inequality with both $L^{p_n}$ norms replaced by $L^p$ holds with loss
\[
 C_{\Sigma,p}[\log(2/\delta)]^{\beta_n(p)},
 \qquad \beta_n(p)=n\left(1-\frac2p\right)\left(\frac5{p_n}+1\right).
\]
Thus the $p=2$ constant is uniform; the weaker $\delta^{-\varepsilon}$ loss follows for every $\varepsilon>0$ throughout $2\le p\le p_n$. Fixed invertible affine images of these surfaces are allowed.
\end{theorem}

The proof is an exact frequency-slicing argument. The slices have a common transverse interval of width comparable to $\delta$, so tensoring with a single one-dimensional function preserves the pointwise square function. The theorem concerns these fixed quadrics. An extension to arbitrary smooth elliptic hypersurfaces would require an additional stability argument, as discussed in Section~\ref{sec:elliptic}.

\subsection{Earlier work and the source of the present argument}
Guth, Wang, and Zhang \cite{GWZ} proved the sharp $L^4$ reverse square function estimate for the cone in $\R^3$, with an arbitrarily small power loss. Cushman, Demeter, and Wu \cite{CDW} obtained the sharp logarithmic power $1/4$ for canonical sectors in that dimension. The powers $5/p_n$ and $\beta_n(p)$ in Theorem~\ref{thm:main} are uniform all-dimensional upper bounds and are not claimed to be optimal. In particular, the known three-dimensional cone estimate is quantitatively stronger.

The four-dimensional cone estimate in \cite{Temur26} was obtained from the scalar local-smoothing and regularity-controlled propagation statements in the September 2026 OpenAI manuscript. The proof here uses the October manuscript \cite{OAI26Oct} in a different way. Besides its scalar critical local-smoothing estimate, we use the Gaussian current, the fiberwise Dirichlet problem, the variational energy, and the energy identity from its Sections 2--5 and Appendix A. The precise statements needed are recorded in Section~\ref{sec:inputs}. We do not use the earlier round-cone propagation theorem. See also \cite{FT}.

This dependence matters. Randomization of scalar local smoothing estimates the square function of free waves in terms of the square function of their data. It does not give a reverse estimate for a coherent sum of angular pieces. The comparison of the coherent and diagonal currents supplies that step. The results of this article are deductions from the stated October inputs; their proofs are not independent proofs of those inputs. The bibliography identifies the immutable October release used here.

\subsection{The proof}
After a bounded angular coloring, the coherent and diagonal data have the same $L^2$ normalization. We construct their currents with the same fiber, the same rank-one background, and one copy of the same positive regularizer. The Dirichlet energy is monotone in the density and has a fixed-background amplification bound.

Put $\ell_h=\log(2/h)$ and $S=A\ell_h$. A dilation to the original fixed source torus changes the wavelength from $h$ to $\widetilde h=h/S$, while the angular pieces retain width $\sqrt h$. At auxiliary radius one, only a bounded number of these pieces are close to the Gaussian center. The remaining sum has exponential decay $e^{-cS}$ and is absorbed by the common regularizer. Thus the coherent energy is bounded by a constant multiple of the diagonal energy on that fiber. Two radial comparisons, with one common time weight, cost $\ell_h^2$ in the $p_n$-th power. This proves the periodic free-wave estimate with norm loss $\ell_h^{2/p_n}$ on time intervals of length $S$.

We then localize to Euclidean space using cutoffs transported along the cap directions. One common anchor time cancels all volume and time-length factors in this localization. Exact temporal reproduction on unit output slabs and vector local smoothing cost
$S^{1/p_n}[\log(2+S/h)]^{2/p_n}$. The total loss is $\ell_h^{5/p_n}$. Positive averaging costs one further logarithm when extending the endpoint operator to unrestricted $L^{p_n}(\ell^2)$. Interpolation there with its uniform $L^2$ bound yields the lower-exponent loss $[\log(2R)]^{\beta_n(p)}$. Finally, exact affine slices give Theorem~\ref{thm:ellipticintro}.

Proposition~\ref{prop:sharpness} shows that every exponent above $p_n$ has a polynomial obstruction for the cone, and that its endpoint constant is at least $c[\log R]^{1/p_n}$. These cone lower bounds are independent of the October inputs. They do not, by themselves, give lower bounds for the elliptic slices.

\section{The cone and its angular decomposition}
\label{sec:target}
Throughout, $d=n+1$ is the ambient spacetime dimension and $n\ge2$ is the spatial dimension. The critical exponent and the associated density exponent are
\begin{equation}\label{eq:exponents}
 \pn=\frac{2n}{n-1}=\frac{2(d-1)}{d-2},\qquad
 \qn=\frac{\pn}{2}=\frac{n}{n-1},\qquad \qn'=n.
\end{equation}
Thus $(d,\pn)=(3,4),(4,3),(5,8/3),(6,5/2)$ in the first four cases. We reserve $q_n$ for the exponent of a squared density.

Let
\[
 \Gamma_n=\{(\xi,\tau)\in\R^n\times\R:
       \tau=|\xi|,\ 1\le|\xi|\le2\}.
\]
For $R\ge2$, set $h=R^{-1}$ and choose a maximal $h^{1/2}$-separated angular net $\{\omega_\theta\}\subset S^{n-1}$. With a fixed enlargement constant $C_0$, put
\[
 \theta=\left\{(\xi,\tau):1\le|\xi|\le2,\quad
 \left|\xi/|\xi|-\omega_\theta\right|\le C_0h^{1/2},\quad
 |\tau-|\xi||\le C_0h\right\}.
\]
There are $O(h^{-(n-1)/2})$ sectors, with bounded overlap. 

\subsection{A smooth partition}
Choose nonnegative smooth bumps $b_\theta$ on $S^{n-1}$, equal to one on the $\sqrt h$-balls about $\omega_\theta$ and supported on their fixed enlargements. The denominator $\sum_\theta b_\theta$ is bounded above and below. Thus
\[
 \psi_\theta=\frac{b_\theta}{\sum_\nu b_\nu},\qquad
 \sum_\theta\psi_\theta=1,\qquad
 |\nabla^k\psi_\theta|\lesssim_k h^{-k/2}
\]
is a smooth partition of unity. If $f$ has Fourier support in the thickened cone, one may set
\[
 \widehat f_\theta(\xi,\tau)
 =\psi_\theta(\xi/|\xi|)\widehat f(\xi,\tau).
\]
The supports form a bounded-overlap cover, rather than a disjoint rectangular partition. The theorem also allows arbitrary functions with these supports, without requiring that they arise from this construction.

The functions are arbitrary within their sectors. They need not be free waves and need not be the pieces of a particular multiplier partition. The right side is the norm of the \emph{pointwise} square function. In particular,
\[
 \left\|\left(\sum_\theta|f_\theta|^2\right)^{1/2}\right\|_p
 \le \left(\sum_\theta\|f_\theta\|_p^2\right)^{1/2},\qquad p\ge2,
\]
so an estimate with the larger right side is not a proof of \eqref{eq:maintheorem}.

\subsection{Frequency planks and dual scales}
Fix $\omega\in S^{n-1}$ and an orthonormal basis $e_1,\ldots,e_{n-1}$ of $\omega^\perp$. The vectors
\[
 g_\omega=2^{-1/2}(\omega,1),\qquad
 E_j=(e_j,0),\qquad n_\omega=2^{-1/2}(\omega,-1)
\]
form an orthonormal frame in $\R^{n+1}$. For $(\xi,\tau)=(r\omega',r+u)$ in the sector about $\omega$, and $q_\omega=(3\omega/2,3/2)$, we have
\[
 |((\xi,\tau)-q_\omega)\cdot g_\omega|\lesssim1,\qquad
 |((\xi,\tau)-q_\omega)\cdot E_j|\lesssim\sqrt h,
\]
and
\[
 |((\xi,\tau)-q_\omega)\cdot n_\omega|
 =2^{-1/2}|r(\omega'\cdot\omega-1)-u|\lesssim h.
\]
Here $1-\omega'\cdot\omega=|\omega'-\omega|^2/2$. Conversely, a sufficiently small plank with these dimensions, on a fixed positive radial band, lies in a fixed sector enlargement. Fixed annular subdivisions are harmless.

\subsection{Fourier convention and normalization}
We use $\widehat f(\xi)=\int e^{-ix\cdot\xi}f(x)\dd x$, suppressing harmless inversion constants. Write $U(t)=e^{it|D|}$; conjugation $U(t)g=\overline{e^{-it|D|}\overline g}$ identifies its norm estimates with the sign convention in the October source, without reversing the time windows. A frequency sector has dimensions
\[
 1\ \times\ (h^{1/2})^{n-1}\ \times\ h.
\]
Indeed, a displacement $h^{1/2}$ tangent to $S^{n-1}$ changes its supporting normal coordinate by $O(h)$, while the radial generator remains flat. Its dual physical plank has lengths $1,(h^{-1/2})^{n-1},h^{-1}$ and volume $h^{-(n+1)/2}$. Under $F(X)=f(RX)$, the frequency becomes $h^{-1}$ and the cone thickness becomes $O(1)$, while the physical packet dimensions become
\begin{equation}\label{eq:packetvolume}
 h\ \times\ (h^{1/2})^{n-1}\ \times\ 1,
 \qquad v_h=h^{(n+1)/2}.
\end{equation}
Both sides of \eqref{eq:maintheorem} acquire the same norm factor $R^{-(n+1)/\pn}$.

\subsection{Disjoint angular groups and shared normalization}
\begin{lemma}[Disjoint angular groups and common normalization]\label{lem:color}
The sector family can be partitioned into $O_{n,C_0}(1)$ groups, each contained in one fixed small source chart and having pairwise disjoint spatial Fourier supports. For such a group of spatial data,
\begin{equation}\label{eq:colornorm}
 \left\|\sum_\theta g_\theta\right\|_2^2=\sum_\theta\|g_\theta\|_2^2.
\end{equation}
A bound proved uniformly for these groups implies the full-family bound with only a dimensional constant.
\end{lemma}
\begin{proof}
Join two angular caps when the closures of fixed enlargements of their spatial projections intersect. Separation of the cap centers bounds every vertex degree by a constant depending only on $n$ and $C_0$, so greedy coloring gives a bounded number of disjoint-support groups. A fixed finite angular cover by smaller source charts has a positive Lebesgue number. For sufficiently small $h$, each complete cap lies in one chart; assign it to one, without applying a new multiplier or splitting its function. Group by this assignment and color. The remaining bounded range of $h$ has only a bounded number of caps and is elementary. Plancherel gives \eqref{eq:colornorm}. Finally use the triangle inequality over the finitely many groups, whose square functions are pointwise bounded by the original one. Thus nonzero data may be normalized to one shared $L^2$ energy without allowing cancellation to separate that energy from the diagonal cap energy.
\end{proof}

\section{Preliminaries}
\label{sec:inputs}
We state the results imported from \cite{OAI26Oct} in the forms used below. The wave equation has $n$ spatial variables. The torus in the auxiliary construction has side $100$, and its measure is ordinary Lebesgue measure. The source construction is used only on this fixed torus and on its original auxiliary domain.

\subsection{Critical local smoothing}
\begin{externalinput}[October critical annular local smoothing]\label{input:ls}
Theorem 1.2 of the October 25, 2026 OpenAI manuscript \cite{OAI26Oct}, printed page 4, states the following. For $n\ge2$, $\pn=2n/(n-1)$, and a fixed smooth annular multiplier $P_N=\phi(D/N)$,
\begin{equation}\label{eq:ls}
 \|U(T)P_Ng\|_{L^{\pn}(\R^n\times[1,2])}
 \lesssim_{n,\phi}[\log(2+N)]^{2/\pn}\|g\|_{L^{\pn}(\R^n)}.
\end{equation}
For $2\le p<\infty$ and $\eps>0$, Theorem 1.1 states, more generally,
\[
 \|U(T)g\|_{L^p(\R^n\times[1,2])}
 \lesssim_{n,p,\eps}\|J^{s_+(p)+\eps}g\|_p,
 \qquad s_+(p)=\max\{0,(n-1)/2-n/p\},
\]
with $\eps=0$ allowed at $p=2$.
\end{externalinput}

\subsection{The auxiliary fiber and the Gaussian current}
The next two inputs summarize the constructions and identities of Sections 2--5 and Appendix A of \cite{OAI26Oct}. The notation $\mathsf H$ refers to the spatial quaternionic Hermitian Hessian in that construction; $H$ is its full Hessian, including the time variable. Their determinant and cofactor are those used in the cited identities. We use the source differential operators only through the following statements.

For a quaternionic matrix, $A^*$ denotes its conjugate transpose. A positive quaternionic current means a matrix field satisfying $A^*=A$ and $\Re(v^*Av)\ge0$ for every quaternionic vector $v$, pointwise.

\begin{externalinput}[Fiber geometry and lifted current]\label{input:current}
Fix $m=n-1$, $a=Kh$, with $K$ the sufficiently large dimensional constant in the source. Let $g(T)=e^{-iT|D|}f$ have annular frequencies comparable to $h^{-1}$, contained in one fixed admissible angular chart, and $\|f\|_2=1$.
Here $T\in\R$ is physical time and $\tau\in\C\simeq\R^2$ is auxiliary time; this $\tau$ is not the Fourier variable in \eqref{eq:sectorintro}. We write $\Delta_{T,\tau}=\partial_T^2+\Delta_{\Re\tau,\Im\tau}$. At fixed $(T,\tau)$, with $s=|\tau|$, the source fiber $\Omega_{T,\tau}$ has real spatial variable $X\in\mathbb T^n$ and auxiliary variables $(P,Z)\in\C^m\times\C$. Put
\[
 Y=(\tau+Z)/\sqrt2,\qquad W=P/Y,\qquad
 \Phi=(Z-\tau)/\sqrt2+P\cdot P/(2Y),\qquad Q=\Phi/a,
\]
where $a=Kh$. The fiber is
\[
 r=|W|^2/w_0^2+100|Q|^2-1<0,\qquad w_0=1/10.
\]
Equivalently, in the fixed ellipsoid coordinates $(W,Q)$,
\[
 Y=\frac{\sqrt2\tau+aQ}{1+W\cdot W/2},\qquad P=WY,
 \qquad Z=\sqrt2Y-\tau.
\]
Spatial quaternionic components are ordered as $(Z,P_1,\ldots,P_m)$; the differential of the displayed $(P,Z)$ chart is used with the corresponding coordinate permutation. Auxiliary complex numbers are embedded in the quaternionic plane $\operatorname{span}_{\R}\{1,\mathbf k\}$.
Source Lemma 4.1 gives, uniformly on the closed fibers for $a<s<4$,
\begin{equation}\label{eq:fiberinterface}
 |Y|\asymp s,\quad
 \dd P\dd Z\asymp a^2s^{2m}\dd W\dd Q,\quad
 \det\mathsf H(r)\asymp a^{-2}s^{-2m},\quad
 \mathsf H(r)\ge cs^{-2}I.
\end{equation}
The full cofactor $\mathsf C(Hr)$ is positive, conserved, tangent to the fiber wall, and has time entry $\det\mathsf H(r)$. The background
\[
 B_0=\lambda|\Phi|^2,\qquad B_\lambda=\mathsf H(B_0),
 \qquad \lambda=h^{-30n}
\]
has spatial rank one and nonzero eigenvalue comparable to $\lambda$.

Source Proposition 2.2 constructs a lifted amplitude $\chi_g$ additive and real-linear in $g$ and a raw current $A^0_g=\chi_g\chi_g^*$ quadratic in $g$. The source lift contains complex conjugation; only additivity is needed below. The raw current is positive, conserved, and annihilates $H(|\Phi|^2)$. Its time density is
\begin{equation}\label{eq:rawtimedensity}
 (A^0_g)_{TT}=|S_0g|^2+|S_1g|^2,
\end{equation}
where each scalar multiplier on a Fourier mode has a bounded polarization factor times
\begin{equation}\label{eq:gaussianmultiplier}
 (\kappa/y)^{m/2}\exp\{-\kappa y|W-v|^2/2\},
 \qquad y=|Y|\asymp s,\quad\kappa\asymp h^{-1}.
\end{equation}
Here $v\in\R^m$ is the mode slope in the fixed chart. The slope map and angular coordinates are uniformly bi-Lipschitz there.

The source smoothing kernel $J_h$ is nonnegative and normalized. If $\vartheta(Q)$ is its fixed passive cutoff, then
\[
 A=\vartheta(Q)J_h*_X A^0_g+h^{50n}\mathsf C(Hr)
\]
is positive, has divergence-free real columns, and satisfies $A:H(B_0)=0$. Its time density is strictly positive. In fiber integrals below, $\dd\zeta=\dd X\dd P\dd Z$ denotes real Lebesgue measure, with periodic measure in $X$. The regularizing cofactor is tangent to the fiber wall, so it contributes no wall flux.
\end{externalinput}

\subsection{The Dirichlet energy}
\begin{externalinput}[Dirichlet problem and energy identity]\label{input:dirichlet}
On a fiber from External input~\ref{input:current}, put $\lambda=h^{-30n}$ and $B_\lambda=\mathsf H(\lambda|\Phi|^2)$.
For any smooth density $\mu$ strictly positive on this closed fiber, source Theorem A.1 supplies a unique smooth zero-boundary potential $u\le0$ satisfying
\[
 D_u=B_\lambda+\mathsf H(u)>0,\qquad \det D_u=\lambda\mu.
\]
Define
\begin{equation}\label{eq:energyfunctional}
 \mathcal F_\lambda(\mu)=\max_v\left\{
 \frac1\lambda\int_{\Omega}v\int_0^1\det(B_\lambda+t\mathsf H(v))\dd t\dd\zeta
                  -\int_{\Omega}v\mu\dd\zeta\right\}.
\end{equation}
The class consists of smooth zero-boundary $v$ with $B_\lambda+\mathsf H(v)>0$; its unique maximizer is $u$. This is source Section 4.4. With $X(\mu)=\int(-u)\mu$, source Lemma 4.3 gives
\begin{equation}\label{eq:energycomparison}
 \frac{n-1}{n+1}X(\mu)\le\mathcal F_\lambda(\mu)
 \le\frac n{n+1}X(\mu).
\end{equation}

We also state the source energy balance explicitly. When $\mu=A_{TT}$ for a positive Hermitian current with divergence-free real columns and $A:H(B_0)=0$, write
\[
 H(B_0+u)=\begin{pmatrix}\ell&c^*\\c&D_u\end{pmatrix},\qquad
 \Sigma=\begin{pmatrix}c^*D_u^{-1}c&c^*\\c&D_u\end{pmatrix}\ge0.
\]
Source Proposition 4.2 gives
\begin{equation}\label{eq:sourceenergybalance}
 \Delta_{T,\tau}\mathcal F+\int_\Omega A:\Sigma\dd\zeta
 =\int \rho_\partial\,G(r)^*A_wG(r)\,\delta(r)\dd\zeta,
 \qquad \rho_\partial=(u/r)|_{r=0}.
\end{equation}
Here $G(r)$ is the full quaternionic gradient of the moving wall, $A_w$ is the current without the cofactor regularizer, and $\delta(r)\dd\zeta=\dd S_{\rm fib}/|\nabla_{\rm fib}r|$. The nonnegative interior contraction makes the wall flux an upper bound for $\Delta\mathcal F$. This identity requires Hermitian symmetry, positivity, conservation, and background annihilation of $A$, not rank one of $A$; the rank-one object is the spatial background. It therefore applies unchanged to the diagonal current.

The constructions are smooth in the displayed parameters. In particular the Dirichlet solution depends smoothly on a smooth one-parameter family of positive densities on a fixed fiber.
\end{externalinput}

\begin{externalinput}[Fine-radius bounds and radial differential inequality]\label{input:fineradius}
For the normalized scalar wave and current above, Proposition 5.1 of \cite{OAI26Oct} gives
\begin{align*}
 \mathcal F(T,\tau)&\gtrsim a^2\|g(T)\|_{p_n}^{p_n} &&(|\tau|=2a),\\
 \mathcal F(T,\tau)&\lesssim a^2\|g(T)\|_{p_n}^{p_n} &&(|\tau|=4a),\\
 \Delta_{T,\tau}\mathcal F(T,\tau)&\le C|\tau|^{-2}e^{-c|\tau|/h}\mathcal F(T,\tau)
 &&(2a\le|\tau|\le3).
\end{align*}
The proof uses scalar multiplier bounds with Gaussian profile \eqref{eq:gaussianmultiplier}, positive-kernel smoothing, comparison barriers for the Dirichlet equation, and the energy identity \eqref{eq:sourceenergybalance}. Those bounds are applied componentwise in Lemma~\ref{lem:hilbertcurrent} below.
\end{externalinput}

\subsection{The spatial maximal theorem of Di Plinio-Parissis}
This final input is used for the lower-exponent range only. It is not used in the endpoint proof.
\begin{externalinput}[Codimension-one plate maximal estimate of \cite{DPP}]\label{input:plates}
Let $0<\delta<1/2$. In $\R^n$, let $\mathcal A_{\omega,\delta}$ be normalized averaging over a plate of thickness $\delta$ in direction $\omega$ and unit width in the orthogonal directions. The codimension-one case of Theorem C in \cite{DPP} gives
\[
 \left\|\sup_{\omega\in S^{n-1}}\mathcal A_{\omega,\delta}|b|\right\|_2
 \lesssim_n\log(2/\delta)\|b\|_2.
\]
Interpolation with $L^\infty$ gives the $L^n$ norm bound $C_n[\log(2/\delta)]^{2/n}$. Fixed positive Schwartz plate majorants satisfy the same bounds by summable common dilations.
\end{externalinput}

\section{Coherent and diagonal energies}
\label{sec:energyproof}
The comparison takes place on one fiber with one fixed background. In particular, no almost-orthogonality statement between different auxiliary radii is assumed.

For this section write $U_-(T)=e^{-iT|D|}$, the source sign convention. Since
$U(T)g=\overline{U_-(T)\overline g}$, conjugation transfers the conclusions to the convention in Section~\ref{sec:coarseroute}, without reversing the time windows. Apply Lemma~\ref{lem:color} first. Thus the data lie in one fixed admissible angular chart, their spatial Fourier supports are disjoint, and we normalize
\begin{equation}\label{eq:sharednormalization}
 \left\|\sum_\theta f_\theta\right\|_2^2
 =\sum_\theta\|f_\theta\|_2^2=1.
\end{equation}
Set $m=n-1$, $g_\theta(T)=U_-(T)f_\theta$, $g=\sum_\theta g_\theta$, and
$S_T=(\sum_\theta|g_\theta(T)|^2)^{1/2}$.

Let $J_h$ be the source's positive normalized spatial smoothing kernel and $0\le\vartheta(Q)\le1$ its fixed passive cutoff. Define
\begin{equation}\label{eq:twocurrents}
 \begin{split}
 A_{\rm coh}&=\vartheta(Q)J_h*_X A^0_g+h^{50n}\mathsf C(Hr),\\
 A_{\rm diag}&=\vartheta(Q)J_h*_X\sum_\theta A^0_{g_\theta}
                           +h^{50n}\mathsf C(Hr).
 \end{split}
\end{equation}
Each has \emph{one} copy of the same regularizer. Set $\mu_{\rm coh}=(A_{\rm coh})_{TT}$ and $\mu_{\rm diag}=(A_{\rm diag})_{TT}$. In particular $\mu>0$. Conservation and background annihilation hold for the diagonal current by linearity in the current. Positivity is preserved by finite sums and averaging with nonnegative weights.

Both currents are Hermitian. Indeed, $(\chi\chi^*)^*=\chi\chi^*$ for every lifted amplitude $\chi$. Since $Hr$ is Hermitian, its cofactor is also Hermitian: for invertible Hermitian $M$, $\mathsf C(M)=(\det M)M^{-1}$, and the singular case follows by continuity. Conjugate transpose commutes with integration against real weights and with finite sums. As $\vartheta(Q)$, $J_h$, and $h^{50n}$ are real, \eqref{eq:twocurrents} gives $A_{\rm coh}^*=A_{\rm coh}$ and $A_{\rm diag}^*=A_{\rm diag}$.

Let $\mathcal F_{\rm coh}=\mathcal F_\lambda(\mu_{\rm coh})$ and $\mathcal F_{\rm diag}=\mathcal F_\lambda(\mu_{\rm diag})$, with the functional in \eqref{eq:energyfunctional}.

\subsection{The source estimates extend to the diagonal current}
\begin{lemma}[Finite-sum current estimates]\label{lem:hilbertcurrent}
With the shared normalization \eqref{eq:sharednormalization}, the diagonal energy satisfies
\begin{align}
 \mathcal F_{\rm diag}(T,\tau)&\lesssim a^2\|S_T\|_{\pn}^{\pn}
                        &&(|\tau|=4a),\label{eq:diaginner}\\
 \Delta_{T,\tau}\mathcal F_{\rm diag}(T,\tau)
 &\le Cs^{-2}e^{-cs/h}\mathcal F_{\rm diag}(T,\tau)
                        &&(2a\le s\le3).\label{eq:diagpotential}
\end{align}
The constants do not depend on the number or individual aperture of the components. Only containment in the fixed source chart and the shared normalization are used. The coherent energy satisfies the source lower bound
$\mathcal F_{\rm coh}(T,\tau)\gtrsim a^2\|g(T)\|_{\pn}^{\pn}$ on $|\tau|=2a$ and the same differential inequality.
\end{lemma}
\begin{proof}
The coherent assertions are Proposition 5.1 of the source. We verify the finite-sum extension of its proof, rather than assuming a vector-valued version of its scalar conclusion.

Regard $G_T=(g_\theta(T))_\theta$ as a finite-dimensional Hilbert-valued wave. Every scalar multiplier in the current acts componentwise, with componentwise conjugation where present in the source lift. Conjugation preserves all Hilbert norms below. If the kernel is bounded by $B J_h$, weighted Cauchy--Schwarz gives
\[
 \|SG_T(X)\|_{\ell^2}^2\lesssim B^2(J_h*\|G_T\|_{\ell^2}^2)(X),
\]
without a factor equal to the Hilbert-space dimension. Put
$P_*=(J_{Ds}*S_T^2)^{1/m}$ with the same fixed large $D$ as in the source. Positive-kernel differentiation and Young give
\[
 |\nabla_X^jP_*|\lesssim(Ds)^{-j}P_*\ (j=1,2),\qquad
 \int P_*^{m+1}\lesssim\int S_T^{\pn}.
\]
Since $\int S_T^2=1$, the source lower bound on the periodized kernel absorbs the single $h^{50n}$ regularizer exactly as in its scalar proof. For $\sigma=s/h$, the multiplier estimates and the additional $J_h$ smoothing therefore give, with fixed finite powers denoted by $M$,
\[
 \mu_{\rm diag}\lesssim s^{-2m}\sigma^M P_*^m,
 \qquad \|A_{{\rm diag},w}\|\lesssim s^{-2m}\sigma^M e^{-c_0\sigma}P_*^m
 \quad\hbox{on the flux wall}.
\]
The normalized vector Fourier coefficients obey
$\sum_\xi\|c_\xi\|_{\ell^2}\lesssim h^{-n/2}$, because there are $O(h^{-n})$ modes. Hence the density ceiling used in source equation (4.7) is unchanged and remains far below the fixed background scale $\lambda$.

The geometric inequality $\mathsf H(rP_*)\ge(1-C/D)P_*\mathsf H(r)$ is unchanged. The same barrier $C\sigma^{M_1}rP_*$ is below the diagonal Dirichlet potential, giving
$0\le u/r\lesssim\sigma^{M_1}P_*$. Combining the density bound, the fiber volume in \eqref{eq:fiberinterface}, and \eqref{eq:energycomparison} gives
\[
 \mathcal F_{\rm diag}(T,\tau)\lesssim a^2\sigma^M\|S_T\|_{\pn}^{\pn}.
\]
At $s=4Kh$ the scaled radius is fixed, proving \eqref{eq:diaginner}. The conserved-current energy balance of source Proposition 4.2 applies to the diagonal current because its Hermitian symmetry, positivity, conservation, and background annihilation were verified after \eqref{eq:twocurrents}. Its wall estimate becomes
\begin{equation}\label{eq:diagflux}
 \Delta\mathcal F_{\rm diag}\lesssim (a/s)^2\sigma^Me^{-c_0\sigma}
                                      \|S_T\|_{\pn}^{\pn}.
\end{equation}

For the matching lower bound, use the source recovery patch $|W|\le\gamma$, $|Q|\le1/400$. Its inverse Gaussian multiplier has kernel constant
$B_{\rm inv}\lesssim_\gamma(hs)^{m/2}\sigma^M e^{C\gamma^2\sigma}$. Apply the inverse multiplier componentwise and sum. On every comparison ball of radius $h$, the diagonal density is therefore at least
\[
 m_X=c_\gamma(hs)^{-m}\sigma^{-M}e^{-C\gamma^2\sigma}S_T(X)^2.
\]
There is no sum of separate $\pn$-th powers at this step: the prescribed density is the full sum of cap squares. The rank-one background leaves $m=n-1$ active Hessian directions: comparison depth scales as $\mu^{1/m}$, energy as $\mu^{1+1/m}$, and $\mu\sim|g|^2$ therefore gives $p=2+2/m=2(d-1)/(d-2)$. The same quadratic comparison on the ball gives
$(-u)\mu_{\rm diag}\gtrsim h^2m_X^{1+1/m}$. Integration on the recovery patch and \eqref{eq:energycomparison} give
\[
 \mathcal F_{\rm diag}(T,\tau)
 \gtrsim_\gamma a^2\sigma^{-M}e^{-C_1\gamma^2\sigma}\|S_T\|_{\pn}^{\pn}.
\]
Choose the fixed aperture $\gamma$ so that $C_1\gamma^2<c_0/2$, as in the source. Dividing \eqref{eq:diagflux} by this lower bound and absorbing the fixed polynomial in the exponential proves \eqref{eq:diagpotential}. The choice of $K$ and all profile constants is the same dimensional choice as in Proposition 5.1.
\end{proof}

\subsection{Monotonicity on a fixed fiber}
\begin{lemma}[Fixed-background variational amplification]\label{lem:energygrowth}
On a fixed source fiber, $\mathcal F_\lambda$ is monotone in the positive density and
\begin{equation}\label{eq:energyamplification}
 \mathcal F_\lambda(L\mu)
 \le L^{(n+1)/(n-1)}\mathcal F_\lambda(\mu),\qquad L\ge1.
\end{equation}
\end{lemma}
\begin{proof}
Every admissible zero-boundary potential $v$ in \eqref{eq:energyfunctional} is nonpositive. At a positive interior maximum its spatial Hessian is nonpositive, whereas $B_\lambda$ has rank one; their sum cannot be positive definite when $n\ge2$. The admissible class is independent of $\mu$, so the variational formula is monotone in $\mu$.

At fixed $\lambda$, let $u_t$ solve the Dirichlet equation with density $t\mu$ and put $F(t)=\mathcal F_\lambda(t\mu)$ and $X(t)=\int(-u_t)t\mu$. Smooth parameter dependence from source Theorem A.1 and stationarity at the maximizer give the envelope identity
$tF'(t)=X(t)$. By \eqref{eq:energycomparison},
\[
 \frac{\dd}{\dd t}\log F(t)
 =\frac{X(t)}{tF(t)}\le\frac{n+1}{n-1}\frac1t.
\]
Integrate from $1$ to $L$ to prove \eqref{eq:energyamplification}. This fixed-fiber step uses arbitrary positive densities, not an assertion that $t\mu$ is the time component of a normalized wave current. The background parameter never changes.
\end{proof}

\subsection{Two radial comparisons at one common time weight}
\begin{lemma}[Elementary radial comparisons]\label{lem:radial}
Let $F\ge0$ satisfy $\Delta_{T,\tau}F\le V(|\tau|)F$ on $b\le|\tau|\le3$, where $b=2Kh$ and $V(s)\le Cs^{-2}e^{-cs/h}$. Let $H$ be its circular mean, and set
\[
 \omega(T)=\sin(\pi T/3),\qquad
 M(s)=\int_0^3\omega(T)H(T,s)\dd T,\qquad D=\log(3/b).
\]
For all sufficiently small $h$,
\begin{equation}\label{eq:radialpair}
 M(b)\lesssim D M(1),\qquad M(1)\lesssim D M(2b).
\end{equation}
\end{lemma}
\begin{proof}
Twice integrating by parts gives
\[
 \int_0^3\omega\,\partial_T^2H\dd T
 =-\frac{\pi^2}{9}M+\frac\pi3\{H(0,s)+H(3,s)\}.
\]
Since the endpoint traces are nonnegative, the differential inequality implies
$M''+s^{-1}M'\le(V+\pi^2/9)M$.

Here is the common radial comparison bound, including its dependence on $h$. With $x=\log(s/b)$, a positive equality solution vanishing at $b$ becomes
$\psi''=P(x)\psi$, where $P=s^2(V+\pi^2/9)$. Normalize $\psi(0)=0$, $\psi'(0)=1$. Then $\psi\ge x$, and the Volterra identity is
\[
 \frac{\psi(x)}x
 =1+\int_0^x\frac{z(x-z)}x P(z)\frac{\psi(z)}z\dd z.
\]
Since $z(x-z)/x\le\min(z,A-z)$ for $x\le A$, Gronwall gives
\begin{equation}\label{eq:radialvolterra}
 \frac{\psi(A)}A\le
 \exp\left(\int_0^A\min(z,A-z)P(z)\dd z\right)\lesssim1.
\end{equation}
For $A=\log(S/b)$ with fixed $S=1$ or $3$, the two contributions obey
\[
 \int_b^S s\log(s/b)V(s)\dd s\lesssim1,\qquad
 \int_b^S s\log(S/s)\dd s\lesssim1.
\]
Reflection $x\mapsto A-x$ preserves the same bound. Hence the equality solution with endpoint values $1,0$ is bounded below by a fixed multiple of the affine function of $\log s$ with those values. Comparing it with the supersolution $M$ on $[b,3]$, using $M(3)\ge0$, gives
$M(1)\gtrsim M(b)\log3/D$. This proves the first inequality in \eqref{eq:radialpair}.

For the second let $w$ solve equality with $w(b)=0$, $w'(b)=b^{-1}$. Then $w>0$ on $(b,3]$, and
\[
 [s(M'w-Mw')]'\le0,\qquad
 b(M'(b)w(b)-M(b)w'(b))=-M(b)\le0.
\]
Thus $M/w$ is nonincreasing, so
$M(1)\le[w(1)/w(2b)]M(2b)$. By \eqref{eq:radialvolterra},
$w(1)\lesssim\log(1/b)$, whereas $w(2b)\ge\log2$. This proves the second inequality. The Wronskian sign is essential: it is an upper comparison for the diagonal energy, not an inward Harnack assertion.
\end{proof}

\section{The periodic estimate on long time intervals}
\label{sec:periodic}
The Gaussian in the auxiliary current resolves angles at the square root of its wavelength. We apply it at a wavelength smaller than the physical cap scale by a logarithmic factor. Only a bounded number of the original caps then lie near any auxiliary center. This removes the dimension-dependent logarithmic factor from the density comparison.

Throughout this section the spatial frequency is comparable to $h^{-1}$, and the angular pieces have aperture $O(\sqrt h)$ and centers separated at scale $\sqrt h$. Put $\ell_h=\log(2/h)$ and $S=A\ell_h$, where $A$ is a sufficiently large fixed dimensional constant. Use ordinary, unnormalized Lebesgue measure on
\[
 \mathbb T_S^n=\R^n/(100S\Z^n).
\]
All constants below are independent of $h$ and $S=A\ell_h$; they may depend on the fixed annulus, cap enlargement, and profile choices.

\subsection{A two-scale use of the fixed source geometry}
\begin{lemma}[The long-window periodic estimate]\label{lem:largeenergy}
For every time center $T_0$, the cap-supported periodic family satisfies
\begin{equation}\label{eq:largeperiodic}
 \left\|\sum_\theta U(T)f_\theta\right\|_{L^{\pn}(\mathbb T_S^n\times[T_0-S,T_0+S])}
 \lesssim\ell_h^{2/\pn}
 \left\|\left(\sum_\theta|U(T)f_\theta|^2\right)^{1/2}\right\|_{L^{\pn}(\mathbb T_S^n\times[T_0-2S,T_0+2S])}.
\end{equation}
\end{lemma}
\begin{proof}
\textbf{Two distinct small parameters.} Set $\widetilde h=h/S$ and make the common dilation
\[
 x=Sx',\qquad T=T_0+St',\qquad
 F_\theta(x')=[U(T_0)f_\theta](Sx').
\]
Then $U(T)f_\theta(Sx')=U(t')F_\theta(x')$. The rescaled spatial torus has side $100$, the frequency is comparable to $\widetilde h^{-1}$, and the angular caps still have width and separation comparable to $\sqrt h$, not $\sqrt{\widetilde h}$. We apply the source energy interface at parameter $\widetilde h$, with its original fixed geometry, and redo its Gaussian comparison for these broader caps. The estimate will follow directly from the source currents at parameter $\widetilde h$.

Apply the fixed chart/color decomposition and shared $L^2$ normalization. Every component remains in one fixed source macrochart. The source interface and Lemma~\ref{lem:hilbertcurrent} use this macrochart condition, additivity, the total $L^2$ norm, and the number of annular Fourier modes; they impose no canonical aperture or individual cap-profile derivative condition. They therefore give the same coherent lower and diagonal upper inner-energy bounds, and the same differential inequality, at frequency $\widetilde h^{-1}$. Use $\widetilde a=K\widetilde h$, $\widetilde b=2\widetilde a$, background $\widetilde\lambda=\widetilde h^{-30n}$, and precisely one shared regulator $\widetilde h^{50n}\mathsf C(Hr)$. All of $h,S,\widetilde h,\widetilde a,\widetilde\lambda$ are fixed during $(T,\tau)$ differentiation. The source period is the original fixed $100$, so no growing-period extension or uniform elliptic-regularity claim is needed.

\textbf{Constant-factor Gaussian comparison.} At the original coarse radius $s=1$, call a cap near if its center is within $C_{\rm near}\sqrt h$ of $\Re W$. There are only $O_n(1)$ near caps. Choose $C_{\rm near}$ larger than the fixed cap diameter constant. Every omitted Fourier mode is at angular distance at least $c\sqrt h$, so its squared Gaussian factor is at most
\[
 \exp(-ch/\widetilde h)=\exp(-cS).
\]
The torus volume is fixed. Parseval and the annular mode count give
\[
 \sum_\xi\|c_\xi\|_{\ell^2}^2=100^{-n},\qquad
 \#\{\xi\}\lesssim\widetilde h^{-n},\qquad
 \sum_\xi\|c_\xi\|_{\ell^2}\lesssim\widetilde h^{-n/2}.
\]
The lifted time-amplitude bound and regulator density consequently satisfy
\[
 |\chi_{{\rm far},T}|^2\lesssim\widetilde h^{-n-m}e^{-cS},\qquad
 \nu=\widetilde h^{50n}\det\mathsf H(r)
       \gtrsim\widetilde h^{50n-2}.
\]
Their ratio is $O(\widetilde h^{-51n-m+2}e^{-cS})$. For each fixed $A$, sufficiently small $h$ gives
$\log(1/\widetilde h)=\log(S/h)\le2\ell_h$. Choose $A$ so large that $cA>2(51n+m-2)+1$. The ratio then tends to zero. Near-cap Cauchy--Schwarz and absorption of the far part by the one shared regulator give
\[
 \mu_{\rm coh}\lesssim\mu_{\rm diag},\qquad
 \mathcal F_{\rm coh}(T,\tau)\lesssim\mathcal F_{\rm diag}(T,\tau)
 \quad(s=1).
\]
The second inequality is Lemma~\ref{lem:energygrowth} with a constant density factor and unchanged background. The near set is used only for a pointwise inequality and is never differentiated as a current.

\textbf{Two radial passages and dilation.} Let $M_{\rm coh}$ and $M_{\rm diag}$ be the circular means of these energies integrated against $\sin(\pi t'/3)$ on $[0,3]$. Lemma~\ref{lem:radial}, the constant-factor comparison at radius one, and Lemma~\ref{lem:hilbertcurrent} give
\begin{align*}
 \widetilde a^2\int_0^3\sin(\pi t'/3)\|U_-(t')\textstyle\sum_\theta F_\theta\|_{p_n}^{p_n}\dd t'
 &\lesssim M_{\rm coh}(\widetilde b)\\
 &\lesssim \log(2/\widetilde h)M_{\rm coh}(1)\\
 &\lesssim \log(2/\widetilde h)M_{\rm diag}(1)\\
 &\lesssim [\log(2/\widetilde h)]^2M_{\rm diag}(2\widetilde b)\\
 &\lesssim \widetilde a^2[\log(2/\widetilde h)]^2
 \int_0^3\sin(\pi t'/3)
 \left\|(\textstyle\sum_\theta|U_-(t')F_\theta|^2)^{1/2}\right\|_{p_n}^{p_n}\dd t'.
\end{align*}
The inner factors $\widetilde a^2$ cancel. The sine is bounded below on $[1,2]$ and above by one on $[0,3]$. Conjugation changes $U_-$ to $U$ without changing either interval. Thus on the original torus the left interval $[1,2]$ and right interval $[0,3]$ obey the free-wave estimate with norm loss
\[
 [\log(2/\widetilde h)]^{2/\pn}\lesssim\ell_h^{2/\pn}.
\]
To obtain left interval $[-1,1]$ and right interval $[-2,2]$, apply this estimate after time translations $-2$ and $-1$. Their left intervals are $[-1,0]$ and $[0,1]$, and both right intervals lie in $[-2,2]$. Sum the $\pn$-th powers. Time translation is implemented by replacing the free-wave data, without a fixed-time $L^{\pn}$ bound. Finally, the dilation back to $(x,T)$ multiplies both spacetime norms by the same $S^{(n+1)/\pn}$. These factors cancel and give \eqref{eq:largeperiodic}. Bounded chart/color sums and homogeneity finish the proof.
\end{proof}

\section{Euclidean localization and temporal reproduction}
\label{sec:coarseroute}
We now pass from periodic free waves to arbitrary functions supported in a thickened cone. We keep the actual pointwise square function throughout the localization argument. The only remaining logarithmic cost comes from applying scalar local smoothing to the time slices in an exact reproduction formula.

Write $\ell_h=\log(2/h)$ and let $1\le L\le C_{\rm win}\ell_h$, with $C_{\rm win}$ fixed. All functions in this section have spatial frequency comparable to $h^{-1}$ and cap aperture $O(\sqrt h)$. Set
\[
 I_L=[-L,L],\qquad J_L=[-2L,2L],\qquad
 N_h=\#\{\theta\}\lesssim h^{-(n-1)/2}.
\]
The positive spatial cap kernel is
\begin{equation}\label{eq:capkernel}
 B_{\theta,h}(y)=c_Mv_h^{-1}
 \left(1+\frac{|y\cdot\omega_\theta|}{h}
 +\frac{|y-(y\cdot\omega_\theta)\omega_\theta|}{\sqrt h}\right)^{-M},
 \qquad v_h=h^{(n+1)/2},
\end{equation}
where $c_M$ makes its integral one. The fixed decay order $M$ is taken sufficiently large whenever the kernel is used.

\subsection{Cutoffs transported along cap directions}
Choose smooth spatial projectors $P_\theta^0$ and $P_\theta$ such that $P_\theta^0=1$ on the input cap and $P_\theta=1$ on a fixed neighborhood of $\supp P_\theta^0$. Their symbols have uniformly bounded derivatives in
\[
 \xi=h^{-1}\varrho\omega_\theta+h^{-1/2}\zeta,
 \qquad\zeta\perp\omega_\theta,
\]
with $\varrho$ in a compact positive interval. The frequency buffer between $\supp P_\theta^0$ and the complement of the plateau of $P_\theta$ is at least $ch^{-1/2}$. Let $\chi$ be smooth, equal to one on $[-8,8]^n$, and supported in $[-10,10]^n$. Put
\[
 \chi_{z,L}(x)=\chi(x/L-z),\qquad z\in\Z^n.
\]

\begin{lemma}[Scaled cutoff transport]\label{lem:egorov}
For cap-supported $u_\theta$ and $|s|\le4L$,
\[
 U(s)P_\theta\chi_{z,L}U(-s)u_\theta(x)
 =\chi_{z,L}(x+s\omega_\theta)u_\theta(x)
   +E_{\theta,z,s,L}u_\theta(x),
\]
where, for every fixed sufficiently large $M$,
\begin{equation}\label{eq:egorovkernel}
 |E_{\theta,z,s,L}(x,y)|
 \le C_M\sqrt h\left\langle\frac{x+s\omega_\theta}{L}-z\right\rangle^{-M}
 B_{\theta,h}(x-y).
\end{equation}
Consequently,
\begin{equation}\label{eq:egorovsquare}
 \left\|\left(\sum_{z,\theta}|E_{\theta,z,s,L}u_\theta|^2\right)^{1/2}\right\|_{p_n}
 \lesssim h^{(n+1)/(4n)}
 \left\|\left(\sum_\theta|u_\theta|^2\right)^{1/2}\right\|_{p_n}.
\end{equation}
The constants are uniform in the allowed $L$ and $s$.
\end{lemma}
\begin{proof}
Write
\[
 R_\theta(\xi,\nu)=|\xi+\nu|-|\xi|-\omega_\theta\cdot\nu.
\]
On the input support, the error has left symbol
\begin{equation}\label{eq:errorsymbol}
 \int e^{i((x+s\omega_\theta)/L-z)\cdot\mu}\widehat\chi(\mu)
 \{e^{isR_\theta(\xi,\mu/L)}P_\theta(\xi+\mu/L)-1\}\dd\mu,
\end{equation}
multiplied by $P_\theta^0(\xi)$. The factor $L^n$ from the Fourier transform of the cutoff cancels the change-of-variable Jacobian $L^{-n}$.

To check all scale factors, put
\[
 a_h(\varrho,\zeta)=\frac{|\zeta|^2}{\sqrt{\varrho^2+h|\zeta|^2}+\varrho},
 \qquad |\xi|=h^{-1}\varrho+a_h(\varrho,\zeta).
\]
For $\nu=\nu_\parallel\omega_\theta+\nu_\perp$, the identity
\begin{equation}\label{eq:residualidentity}
 R_\theta(\xi,\nu)=a_h(\varrho+h\nu_\parallel,\zeta+\sqrt h\nu_\perp)
                   -a_h(\varrho,\zeta)
\end{equation}
is exact. When $|\mu|\le cLh^{-1/2}$, the two arguments lie in one fixed compact set with positive radial coordinate, and the projector in \eqref{eq:errorsymbol} is one. Uniform smooth bounds for $a_h$ and $|s|/L\le4$ imply
\[
 |sR_\theta(\xi,\mu/L)|
 \lesssim\sqrt h|\mu|+(h/L)|\mu|^2.
\]
Every fixed derivative in $(\varrho,\zeta,\mu)$ of
$e^{isR_\theta(\xi,\mu/L)}-1$ is bounded by $C\sqrt h\langle\mu\rangle^C$; derivatives in $\mu$ gain inverse factors of $L$. On the complementary region, all fixed derivatives of the remaining factors cost at most a polynomial in $L$ and $h^{-1}$, whereas the Schwartz tail starts at $cLh^{-1/2}$. Taking enough decay seminorms of $\widehat\chi$ makes this contribution $O(h^A)$ for any prescribed $A$, uniformly for $L\ge1$. Near $\xi+\mu/L=0$, the product containing $P_\theta(\xi+\mu/L)$ vanishes, so the apparent nonsmoothness of the Euclidean norm is harmless.

Integrating by parts in $\mu$ now gives the scaled spatial weight in \eqref{eq:egorovkernel}, with uniform normalized frequency derivatives. Anisotropic integration by parts in $\xi$ gives the cap kernel. Weighted Cauchy--Schwarz and the bounded lattice sum yield
\[
 \sum_z|E_{\theta,z,s,L}u_\theta(x)|^2
 \lesssim h(B_{\theta,h}*|u_\theta|^2)(x).
\]
To sum in $\theta$, let $q_n=p_n/2$, so $q_n'=n$. Young's inequality and the finite-family bound give
\[
 \|\max_\theta B_{\theta,h}*|b|\|_n
 \le\left(\sum_\theta\|B_{\theta,h}*|b|\|_n^n\right)^{1/n}
 \le N_h^{1/n}\|b\|_n.
\]
Duality therefore bounds the $L^{q_n}$ norm of
$\sum_\theta B_{\theta,h}*|u_\theta|^2$ by
$N_h^{1/n}\|\sum_\theta|u_\theta|^2\|_{q_n}$. Taking square roots gives
\[
 \sqrt h\,N_h^{1/(2n)}\lesssim h^{(n+1)/(4n)},
\]
which proves \eqref{eq:egorovsquare}. No directional maximal estimate is used here.
\end{proof}

\subsection{From the periodic to the Euclidean free-wave estimate}
\begin{lemma}[Euclidean localization]\label{lem:euclidean}
Suppose the periodic inequality
\begin{equation}\label{eq:periodicabstract}
 \left\|\sum_\theta U(t)v_\theta\right\|_{L^{p_n}(\mathbb T_L^n\times I_L)}
 \le A(h,L)
 \left\|\left(\sum_\theta|U(t)v_\theta|^2\right)^{1/2}\right\|_{L^{p_n}(\mathbb T_L^n\times J_L)}
\end{equation}
holds for fixed enlargements of the caps, where $A(h,L)\ge1$. Then
\begin{equation}\label{eq:freewaveSF}
 \left\|\sum_\theta U(t)v_\theta\right\|_{L^{p_n}(\R^n\times I_L)}
 \lesssim A(h,L)
 \left\|\left(\sum_\theta|U(t)v_\theta|^2\right)^{1/2}\right\|_{L^{p_n}(\R^n\times J_L)}.
\end{equation}
The same assertion holds for common translates of both time intervals.
\end{lemma}
\begin{proof}
Let $u_\theta(t)=U(t)v_\theta$ and $S_u=(\sum_\theta|u_\theta|^2)^{1/2}$. Choose one time $t_a\in J_L$ such that
\begin{equation}\label{eq:anchortime}
 \|S_u(t_a)\|_{p_n}\le(4L)^{-1/p_n}
 \|S_u\|_{L^{p_n}(\R^n\times J_L)}.
\end{equation}
Tile space by the disjoint cubes $Q_{z,L}=L(z+[-1/2,1/2)^n)$ and let
$\mathcal D_{z,L}=Lz+[-50L,50L)^n$. Periodize
\[
 w_{\theta,z}=P_\theta(\chi_{z,L}u_\theta(t_a))
\]
on $\mathcal D_{z,L}$. Applying the projector after the cutoff ensures support in a fixed enlarged cap. The estimate \eqref{eq:periodicabstract} applies to the periodized waves on the same absolute time intervals by taking time-zero data $U_{\mathbb T}(-t_a)w_{\theta,z}$. This step uses no fixed-time $L^{p_n}$ bound for the wave group.

We first control the coherent error on $Q_{z,L}\times I_L$. If $K_{\theta,s}$ is the actual kernel of $U(s)P_\theta$, integration by parts in normalized cap coordinates gives, for every integer $N$,
\[
 |K_{\theta,s}(y)|\lesssim_N h^{-(n+1)/2}(1+|s|)^N
 \left(1+\frac{|y\cdot\omega_\theta+s|}{h}
          +\frac{|y-(y\cdot\omega_\theta)\omega_\theta|}{\sqrt h}\right)^{-N}.
\]
For $|s|\le4L$ and $|y+s\omega_\theta|\ge cL$, it follows that, for every prescribed $M$ and fixed sufficiently large $K$,
\begin{equation}\label{eq:largestrongtails}
 |K_{\theta,s}(y)|\le C_{M,K}h^ML^{-n}\langle y/L\rangle^{-K}.
\end{equation}
Indeed the preceding estimate supplies $h^{N/2-(n+1)/2}$ times a profile at scale $L$. Choose $N$ large and absorb the remaining $L^n$ using $L\lesssim\ell_h$. The omitted cutoff region and all nonprincipal periods stay a fixed multiple of $L$ from the propagation center. The kernel in \eqref{eq:largestrongtails} has integral $O(h^M)$; the sum over nonzero period indices is bounded independently of $L$.

Young's inequality and $\sum_\theta|u_\theta(t_a)|\le N_h^{1/2}S_u(t_a)$ show that the spatial $\ell^{p_n}(z)$ sum of the coherent error norms is at most
$C_Mh^MN_h^C\|S_u(t_a)\|_{p_n}$. The order $M$ may be chosen after the polynomial cap count. Time integration contributes $(2L)^{1/p_n}$, canceled by \eqref{eq:anchortime}. We use the strong tails here, rather than summing the weaker cutoff-transport errors coherently.

For the zero period on the right side, Lemma~\ref{lem:egorov} gives
\[
 U(t-t_a)w_{\theta,z}(x)
 =\chi_{z,L}(x+(t-t_a)\omega_\theta)u_\theta(x,t)
  +E_{\theta,z,t-t_a,L}u_\theta(x,t).
\]
Since $p_n\ge2$ and the cutoffs have bounded overlap, pointwise
\begin{align*}
 \sum_z\left(\sum_\theta
 |\chi_{z,L}(x+(t-t_a)\omega_\theta)u_\theta(x,t)|^2\right)^{p_n/2}
 &\le\left(\sum_{z,\theta}
 |\chi_{z,L}(x+(t-t_a)\omega_\theta)u_\theta(x,t)|^2\right)^{p_n/2}\\
 &\lesssim S_u(x,t)^{p_n}.
\end{align*}
The analogous inequality followed by \eqref{eq:egorovsquare} sums the errors in $\ell^{p_n}(z)$. Integration over $J_L$ involves the actual square function at each time and adds no factor of $L$.

For completeness, the nonprincipal-period volume factors cancel as follows. Put
\[
 E_z=\int\chi_{z,L}(y)^2S_u(y,t_a)^2\dd y.
\]
Kernel Cauchy--Schwarz and \eqref{eq:largestrongtails} bound their pointwise vector norm on $\mathcal D_{z,L}$ by
\[
 C_Mh^ML^{-n/2}E_z^{1/2}.
\]
Jensen and bounded overlap give
\[
 \sum_zE_z^{p_n/2}\lesssim
 L^{n(p_n/2-1)}\|S_u(t_a)\|_{p_n}^{p_n}.
\]
Integration over a period cube contributes $L^{n/p_n}$, which cancels the factors $L^{-n/2}$ and $L^{n/2-n/p_n}$ above. The remaining time factor $(4L)^{1/p_n}$ is canceled by the same anchor. Thus the nonprincipal-period terms are bounded by an arbitrarily small power of $h$ times the desired right side.

Apply \eqref{eq:periodicabstract}, sum in $\ell^{p_n}(z)$, and use the preceding bounds. The total factor is
$CA(h,L)(1+h^{(n+1)/(4n)}+h^M)$, proving \eqref{eq:freewaveSF}. Common time translations are obtained by changing the initial data, not by estimating the wave group at a fixed time.
\end{proof}

\subsection{Vector local smoothing}
\begin{lemma}[Bounded and growing time intervals]\label{lem:vectorls}
For every finite family of annular spatial data and $D\ge1$,
\begin{equation}\label{eq:vectorls}
 \left\|\{U(t)a_\theta\}_\theta\right\|_{L^{\pn}(\R^n\times[-D,D];\ell^2)}
 \lesssim D^{1/\pn}[\log(2+D/h)]^{2/\pn}
 \|\{a_\theta\}_\theta\|_{L^{\pn}(\R^n;\ell^2)}.
\end{equation}
\end{lemma}
\begin{proof}
We first extend External input~\ref{input:ls} from $[1,2]$ to $[-2,2]$ by dilation, not by replacing the data by a time translate. On $[a,2a]$, $2h\le a\le1$, set $x=ay$, $t=as$, and $f_a(y)=f(ay)$. The new frequency is $a/h$, so
\[
 \int_a^{2a}\|U(t)f\|_{\pn}^{\pn}\dd t
 \lesssim a[\log(2+a/h)]^2\|f\|_{\pn}^{\pn}.
\]
The dyadic sum of these factors is $O(\ell_h^2)$. On $|t|\lesssim h$ the rescaled annular multiplier kernel has uniformly bounded $L^1$ norm, and negative times follow by conjugation. Randomize the annular family and apply Khintchine on both sides of this scalar estimate; the constants are independent of the number of components. This proves the vector estimate on a fixed bounded interval. Finally the common dilation $(x,t)=D(y,s)$ changes the frequency to $D/h$ and contributes the norm factor $D^{1/\pn}$, proving \eqref{eq:vectorls}.
\end{proof}

\subsection{Exact reproduction of the normal-frequency band}
\begin{lemma}[Temporal transfer]\label{lem:continuous}\label{lem:largetransfer}
Suppose \eqref{eq:freewaveSF} holds. For arbitrary Schwartz functions $g_\theta$ supported in the high-frequency cone sectors, with normal thickness $O(1)$,
\begin{equation}\label{eq:largetransfer}
 \left\|\sum_\theta g_\theta\right\|_{p_n}
 \lesssim A(h,L)L^{1/p_n}[\log(2+L/h)]^{2/p_n}
 \left\|\left(\sum_\theta|g_\theta|^2\right)^{1/2}\right\|_{p_n}.
\end{equation}
\end{lemma}
\begin{proof}
The spatial Fourier transform of $v_\theta(t)=U(-t)g_\theta(t)$ is obtained by multiplying by $e^{-it|\xi|}$. Its temporal Fourier support is therefore in a fixed interval, since the new temporal frequency is $\tau-|\xi|$. Choose one Schwartz function $\eta$ whose Fourier transform is one on this interval. Temporal convolution gives the exact identity
\begin{equation}\label{eq:continuousreproduce}
 g_\theta(t)=\int_\R\eta(t-s)U(t-s)g_\theta(s)\dd s.
\end{equation}
Thus the data used in the transfer are actual time slices of $g_\theta$; cancellation within a cap and within its normal-frequency band has already taken place.

On a unit output slab $I_j=[j-1/2,j+1/2)$, Minkowski, followed by the free-wave estimate on $[j-L,j+L]$, gives
\begin{align*}
 \left\|\sum_\theta g_\theta\right\|_{L^{p_n}(\R^n\times I_j)}
 \lesssim A(h,L)\int_\R\sup_{t\in I_j}|\eta(t-s)|
 \left\|\{U(t-s)g_\theta(s)\}_\theta\right\|_{L^{p_n}(\R^n\times[j-2L,j+2L];\ell^2)}\dd s.
\end{align*}
The supremum of the reconstruction weight is taken on the unit interval, before enlarging the free-wave norm. Hence
$\sup_{t\in I_j}|\eta(t-s)|\lesssim_N\langle j-s\rangle^{-N}$ independently of $L$.

On the right time window, $|t-s|\le2L+|j-s|$. Lemma~\ref{lem:vectorls} bounds the vector free-wave norm by the time-slice norm times
\[
 (2L+|j-s|)^{1/p_n}
 [\log(2+(2L+|j-s|)/h)]^{2/p_n}
 \lesssim L^{1/p_n}[\log(2+L/h)]^{2/p_n}\langle j-s\rangle^C.
\]
The Schwartz factor absorbs this polynomial. If
$a(s)=\|\{g_\theta(s)\}_\theta\|_{L^{p_n}_x(\ell^2)}$, the slab norm is at most
\[
 CA(h,L)L^{1/p_n}[\log(2+L/h)]^{2/p_n}
 \int_\R\langle j-s\rangle^{-2}a(s)\dd s.
\]
The integral operator maps $L^{p_n}(\dd s)$ to $\ell^{p_n}(j)$ because
\[
 \sup_j\int_\R\langle j-s\rangle^{-2}\dd s
 +\sup_s\sum_j\langle j-s\rangle^{-2}\lesssim1.
\]
Jensen's inequality and summation now prove \eqref{eq:largetransfer}. The output slabs have length one, so no additional overlap factor arises from the long time windows.
\end{proof}

\begin{proof}[Proof of Theorem~\ref{thm:main} at the endpoint]
Set $L=S=A\ell_h$ as in Lemma~\ref{lem:largeenergy}. That lemma and Lemma~\ref{lem:euclidean} give the Euclidean free-wave estimate with
$A(h,L)\lesssim\ell_h^{2/p_n}$. Lemma~\ref{lem:continuous} contributes the remaining factors
\[
 \ell_h^{2/p_n}\,S^{1/p_n}\,[\log(2+S/h)]^{2/p_n}
 \lesssim\ell_h^{5/p_n}.
\]
The chart and color decompositions cost only a dimensional constant. The common physical dilation described in Section~\ref{sec:target} restores the original frequency-one cone and changes both norms by the same factor. This proves the endpoint assertion for arbitrary sector-supported pieces. The lower exponents are proved in the next section.
\end{proof}

\section{Hilbert-valued estimates and the full exponent range}
\label{sec:corollaries}
For the lower-exponent interpolation argument we use External input~\ref{input:plates}. This is the only part of the cone proof that requires the spatial maximal theorem.

\subsection{A dimension-free Hilbert-valued consequence}
\begin{corollary}[Hilbert-valued endpoint]\label{cor:hilbertendpoint}
Let $H$ be a Hilbert space. For sector-supported $H$-valued Schwartz functions,
\[
 \left\|\sum_\theta f_\theta\right\|_{L^{\pn}(\R^d;H)}
 \lesssim[\log(2R)]^{5/p_n}
 \left\|\left(\sum_\theta\|f_\theta\|_H^2\right)^{1/2}\right\|_{\pn},
\]
with constants independent of $\dim H$. 
\end{corollary}
\begin{proof}
For finitely many channels $f_{\theta,j}$, apply the scalar endpoint estimate to
$F_\theta=\sum_j\varepsilon_jf_{\theta,j}$ with independent random signs. Take the $\pn$-th power and expectation. Scalar Khintchine on the left and Hilbert-valued Khintchine on the right, with Hilbert space $\ell^2_\theta$, give
\[
 \left\|\left(\sum_j\left|\sum_\theta f_{\theta,j}\right|^2\right)^{1/2}\right\|_{\pn}
 \lesssim[\log(2R)]^{5/p_n}
 \left\|\left(\sum_{\theta,j}|f_{\theta,j}|^2\right)^{1/2}\right\|_{\pn}.
\]
Finite-dimensional orthogonal projections and norm convergence give the general Hilbert-space assertion. This randomizes the already-proved reverse estimate; it does not reverse a randomized scalar local-smoothing inequality.
\end{proof}

\subsection{Positive averaging for the auxiliary corollaries}
Using the spatial kernel from Section~\ref{sec:coarseroute}, put
\[
 K_{\theta,h}(y,s)=B_{\theta,h}(y+s\omega_\theta)\rho(s),
 \qquad \rho(s)=c_M(1+|s|)^{-M}.
\]
All kernels are nonnegative and normalized; the fixed decay order is chosen sufficiently large.
\begin{proposition}[Positive averaging]\label{prop:positiveaverage}
Assume External inputs~\ref{input:ls} and~\ref{input:plates}. For all nonnegative families $G_\theta$,
\begin{equation}\label{eq:positiveaverage}
 \left\|\sum_\theta K_{\theta,h}*G_\theta\right\|_{\qn}
 \lesssim_n[\log(2/h)]^2\left\|\sum_\theta G_\theta\right\|_{\qn}.
\end{equation}
Equivalently, taking square roots costs at most one logarithm at exponent $\pn$.
\end{proposition}
\begin{proof}
Write $L_h=\log(2/h)$. External input~\ref{input:plates} gives spatial plate maximal strong $L^2$ norm $O(L_h)$. Interpolation with $L^\infty$ gives $L^n$ norm $O(L_h^{2/n})$. Indeed, after the common spatial dilation $x=\sqrt h X$, the plate of longitudinal width $h$ and transverse width $\sqrt h$ becomes a unit codimension-one plate of thickness $\sqrt h$. Fixed positive Schwartz majorants obey the same maximal bound by summable common spatial dilations.

Lemma~\ref{lem:vectorls} supplies the bounded-time scalar local-smoothing norm loss $L_h^{2/\pn}$.

We need packets with an actual lower bound, not only localization. Let $\lambda_0>0$ be fixed. Choose a nonzero nonnegative smooth Fourier bump $b(\varrho,\zeta)$ supported in $|\varrho|+|\zeta|\le c$, with $c$ fixed sufficiently small, and normalize its inverse transform $\psi$ in $L^2$. In coordinates parallel and perpendicular to $\omega_\theta$, put
\[
 \phi_{\theta,h}(y)=v_h^{-1/2}e^{i\lambda_0y_\parallel/h}
                       \psi(y_\parallel/h,y_\perp/\sqrt h).
\]
Then $\|\phi_{\theta,h}\|_2=1$, and its Fourier support lies in a fixed annulus at frequency $h^{-1}$ and a cap of aperture $O(\sqrt h)$. At $x=-t\omega_\theta+y$, remove the carrier phase. The remaining Fourier phase is
\[
 (y_\parallel/h)\varrho+(y_\perp/\sqrt h)\cdot\zeta
 +t\frac{|\zeta|^2}{\sqrt{(\lambda_0+\varrho)^2+h|\zeta|^2}
                         +\lambda_0+\varrho}.
\]
On $|y_\parallel|\le h$, $|y_\perp|\le\sqrt h$, $|t|\le1$, its modulus is at most $Cc+Cc^2$. Choose $c$ so this is less than $\pi/3$. Since $b\ge0$, the integral has a uniform positive real part. If $A_{\theta,h}$ is the normalized indicator of this spatial plate and $W_{\theta,h}=|\phi_{\theta,h}|^2$, we have
\begin{equation}\label{eq:packetlowerupper}
 |U(t)\phi_{\theta,h}(x)|^2\gtrsim A_{\theta,h}(x+t\omega_\theta),
 \qquad
 W_{\theta,h}(y)\lesssim_M v_h^{-1}
 (1+|y_\parallel|/h+|y_\perp|/\sqrt h)^{-M}.
\end{equation}
All constants are independent of direction, $h$, and $t$ in the stated interval.

Randomize each member of a finite family of translated packets separately, with coefficients $\sqrt{c_{\theta,j}}$. Khintchine and the bounded-time scalar estimate, followed by squaring the norm, give for finite positive atomic measures $\mu_\theta$
\[
 \left\|\sum_\theta\int|U(t)\phi_{\theta,h}(x-y)|^2\dd\mu_\theta(y)
       \right\|_{L^{\qn}_{x,t}([-1,1])}
 \lesssim L_h^{4/\pn}
       \left\|\sum_\theta W_{\theta,h}*\mu_\theta\right\|_{L^{\qn}_x}.
\]
At this stage both sides contain \emph{convolved} atomic measures. Approximate nonnegative compactly supported smooth densities $a_\theta$ by positive Riemann-sum measures. The packet kernels are continuous with uniform spatial decay on bounded time intervals, so the convolved expressions converge in the displayed norms for each fixed $h$. We may therefore replace $\dd\mu_\theta$ by $a_\theta\dd y$ before invoking any $L^{\qn}$ duality on the input. Truncation and density then give arbitrary nonnegative $L^{\qn}$ data.

Since $\qn'=n$, the spatial maximal bound now gives, for these density functions,
\[
 \left\|\sum_\theta W_{\theta,h}*a_\theta\right\|_{\qn}
 \lesssim L_h^{2/n}\left\|\sum_\theta a_\theta\right\|_{\qn}.
\]
Explicitly, pairing against $b\ge0$ in $L^n$ moves each kernel onto $b$ and bounds the pairing by
$\int(\sum_\theta a_\theta)\sup_\theta(\widetilde W_{\theta,h}*b)$.
Use \eqref{eq:packetlowerupper} on the evolved side. Since $4/\pn+2/n=2$, this proves
\begin{equation}\label{eq:positivetransport}
 \left\|\sum_\theta(A_{\theta,h}*a_\theta)(x+t\omega_\theta)
 \right\|_{L^{\qn}(\R^n\times[-1,1])}
 \lesssim L_h^2\left\|\sum_\theta a_\theta\right\|_{L^{\qn}(\R^n)}.
\end{equation}

For spacetime convolution by the normalized hard planks, restrict output time to $I_j=[j,j+1]$; source times then lie in $[j-1,j+2]$. Minkowski, \eqref{eq:positivetransport}, and H\"older on that fixed interval give
\[
 \left\|\sum_\theta K^0_{\theta,h}*G_\theta\right\|_{L^{\qn}(\R^n\times I_j)}
 \lesssim L_h^2
 \left(\int_{j-1}^{j+2}\left\|\sum_\theta G_\theta(\cdot,s)\right\|_{\qn}^{\qn}\dd s\right)^{1/\qn}.
\]
Raise to $\qn$ and sum in $j$ using bounded overlap. This is a local-slab argument, not a global Minkowski bound with an unwanted $L^1$ time norm.

Finally, the kernels $K_{\theta,h}$ in the statement are bounded by a sum of common isotropic dilations of $K^0_{\theta,h}$ with summable coefficients $O(2^{-j(M-n-1)})$, for $M>n+1$. Both the anisotropic widths and the shear $y+s\omega_\theta$ scale compatibly. Each dilated convolution family has the same $L^{\qn}$ operator norm, by dilation of all spatial and temporal coordinates. Summing the coefficients proves \eqref{eq:positiveaverage}.
\end{proof}
For $n>2$, the weak-$L^2$ bound $O(\sqrt{\log(2/h)})$ in the same external theorem gives the sharper square-root loss $\log(2/h)^{1-1/(2n)}$. The uniform one-logarithm form above also covers $n=2$ and suffices here.

\subsection{Interpolation on unrestricted vector spaces}
\begin{corollary}[Full reverse square-function range]\label{cor:range}
For every $2\le p<\pn$, the canonical cone sectors satisfy
\[
 \left\|\sum_\theta f_\theta\right\|_p
 \lesssim_{n,p}[\log(2R)]^{\beta_n(p)}
 \left\|\left(\sum_\theta|f_\theta|^2\right)^{1/2}\right\|_p,
\]
where
\[
 \alpha=\frac{1/2-1/p}{1/2-1/p_n}=n\left(1-\frac2p\right),
 \qquad \beta_n(p)=\alpha\left(\frac5{p_n}+1\right).
\]
At $p=2$ the constant is independent of $R$. At $p=\pn$ use the better direct exponent $5/p_n$ from Theorem~\ref{thm:main}, rather than the interpolated exponent $5/p_n+1$. Consequently, an $R^\eps$ loss is also available for every $\eps>0$ throughout $2\le p\le p_n$. No optimality of these logarithmic powers is asserted.
\end{corollary}
\begin{proof}
One must not interpolate Fourier-supported subspaces without proving the required interpolation identity. Instead work on the unrestricted spaces $L^p(\ell^2)$.

In the high-frequency normalization, choose a smooth spacetime Fourier multiplier $P_\theta$ equal to one on the entire given sector and supported in a fixed enlarged sector. In coordinates
\[
 \xi=h^{-1}\varrho\omega_\theta+h^{-1/2}\zeta,
 \quad\zeta\perp\omega_\theta,\qquad
 \sigma=\tau-\xi\cdot\omega_\theta,
\]
the symbol has uniformly bounded smooth seminorms and fixed compact support. Indeed $\varrho$ stays separated from zero and
\[
 |\xi|-\xi\cdot\omega_\theta
 =\frac{|\zeta|^2}{\sqrt{\varrho^2+h|\zeta|^2}+\varrho}
\]
has uniform derivatives. Since the Fourier pairing is
$(x+t\omega_\theta)\cdot\xi+t\sigma$, integration by parts gives a kernel bound
\[
 |k_{\theta,h}(x,t)|\lesssim_M v_h^{-1}
 \left(1+|t|+\frac{|(x+t\omega_\theta)\cdot\omega_\theta|}{h}
 +\frac{|x-(x\cdot\omega_\theta)\omega_\theta|}{\sqrt h}\right)^{-M}.
\]
Choosing the decay order large enough yields $|k_{\theta,h}|\lesssim K_{\theta,h}$ and $\|k_{\theta,h}\|_1\lesssim1$, for the positive cone-plank majorants already used. Hence, for arbitrary vector inputs $g=(g_\theta)$,
\[
 |P_\theta g_\theta|^2\lesssim K_{\theta,h}*|g_\theta|^2.
\]
Proposition~\ref{prop:positiveaverage} gives
\[
 \|(P_\theta g_\theta)_\theta\|_{L^{\pn}(\ell^2)}
 \lesssim\log(2/h)\|g\|_{L^{\pn}(\ell^2)}.
\]
Define $T_hg=\sum_\theta P_\theta g_\theta$. Apply Theorem~\ref{thm:main} to the supported family $(P_\theta g_\theta)$ to obtain
\[
 \|T_hg\|_{\pn}\lesssim[\log(2/h)]^{5/p_n+1}
                              \|g\|_{L^{\pn}(\ell^2)}.
\]
At $p=2$, Plancherel and bounded overlap of the enlarged supports give
$\|T_hg\|_2^2\lesssim\sum_\theta\|g_\theta\|_2^2$. Vector-valued Riesz--Thorin on the full spaces $L^p(\ell^2)$, with $1/p=(1-\alpha)/2+\alpha/p_n$, now gives the norm loss $[\log(2/h)]^{\alpha(5/p_n+1)}$. Since $1/2-1/p_n=1/(2n)$, this is exactly the stated $\beta_n(p)$. For the original supported family, $P_\theta f_\theta=f_\theta$ exactly. Substitute this family and return to the original scale.
\end{proof}

\section{Square functions for elliptic quadrics}
\label{sec:elliptic}
We pass from the cone in $\R^{n+1}$ to hypersurfaces in $\R^n$.
The exponent is therefore
\[
 p_n=\frac{2n}{n-1},\qquad n\ge2.
\]
The passage is by an exact frequency slice. It preserves the pointwise
square function, as well as the logarithmic loss.

Let $\Sigma\subset\R^n$ be a fixed smooth hypersurface, restricted to a
compact patch when necessary. Put $\delta=R^{-1}$, and choose a
$\sqrt\delta$-separated maximal net $\{z_\theta\}$ on $\Sigma$. A canonical cap is a
fixed enlargement of
\begin{equation}\label{eq:ellipticcapmodel}
 \vartheta_\theta=
 \{z+t\nu_\Sigma(z):z\in\Sigma,\ |z-z_\theta|\le C_0\sqrt\delta,
                  \ |t|\le C_0\delta\},
\end{equation}
where $\nu_\Sigma$ is a unit normal. On a graph
$\Sigma=\{(u,\phi(u)):u\in K\}$, with $K$ compact and $\nabla\phi$ bounded,
these caps are, up to fixed enlargements, the sets
\begin{equation}\label{eq:ellipticgraphcaps}
 |u-u_\theta|\le C_0\sqrt\delta,
 \qquad |v-\phi(u)|\le C_0\delta.
\end{equation}
They have tangential widths comparable to $\sqrt\delta$ and normal width
comparable to $\delta$. Their overlap is bounded independently of
$\delta$. For a patch with boundary we use the caps from a fixed smooth
extension of the patch. All constants below may depend on this extension,
the compact patch, and the fixed enlargement.

We prove Theorem~\ref{thm:ellipticintro} with loss
$[\log(2/\delta)]^{\beta_n(p)}$ for $2\le p<p_n$, where
$\beta_n(p)=n(1-2/p)(5/p_n+1)$, and with the better direct loss
$[\log(2/\delta)]^{5/p_n}$ at $p=p_n$.
The endpoint has precisely the external dependencies of the cone endpoint.
For $2<p<p_n$, the assertion additionally uses the positive-averaging
estimate in Section~\ref{sec:corollaries}; the $p=2$ bound is uniform.

\begin{lemma}[Tensor-product slicing]\label{lem:elliptictensorslice}
Suppose the arbitrary-piece reverse square-function estimate at exponent
$p$ is valid for a family $\{\Theta_\beta\}$ in $\R^{n+1}$ with constant
$A(R)$. Suppose there are an invertible linear map $B$ on $\R^{n+1}$, a
vector $b$, and $c>0$ such that
\[
 b+B(\vartheta_\theta\times[-c\delta,c\delta])
       \subset\Theta_{\beta(\theta)}.
\]
If each index $\beta$ has at most $M$ preimages, then
\[
 \left\|\sum_\theta f_\theta\right\|_{L^p(\R^n)}
 \lesssim_M A(R)
 \left\|\left(\sum_\theta|f_\theta|^2\right)^{1/2}\right\|_{L^p(\R^n)}
\]
for every family supported in $\{\vartheta_\theta\}$.
\end{lemma}
\begin{proof}
Choose a nonzero Schwartz function $\psi$ with
$\operatorname{supp}\widehat\psi\subset[-c,c]$. Set
\[
 G_\theta(x,s)=f_\theta(x)\psi(\delta s),\qquad
 F_\theta(X)=e^{ib\cdot X}G_\theta(B^T X).
\]
Then
\[
 \operatorname{supp}\widehat F_\theta
 \subset b+B(\vartheta_\theta\times[-c\delta,c\delta]).
\]
Group the $F_\theta$ with the same image index. Pointwise,
\[
 \sum_\beta\left|\sum_{\beta(\theta)=\beta}F_\theta\right|^2
 \le M\sum_\theta|F_\theta|^2.
\]
Apply the assumed estimate to the grouped family. The $L^p$ norms of the
sum and the square function both have the factor
\[
 |\det B|^{-1/p}\delta^{-1/p}\|\psi\|_p.
\]
Cancel it. This proves the lemma without a limiting restriction argument
and without applying a multiplier to the functions $f_\theta$.
\end{proof}

We use the cone theorem on any fixed radial band $a\le|\xi|\le b$, where
$0<a<b<\infty$, with fixed sector enlargements. This version follows from
the stated version by finitely many fixed radial subdivisions and common
dilations. To see that the subdivisions do not alter the right side,
choose smooth annular cutoffs $\chi_j(\xi)$ summing to one on the band.
Each $\chi_j(D_x)$ acts on every sector by the same convolution kernel,
whose $L^1$ norm is bounded independently of $R$. Minkowski's inequality
in $\ell^2$ and Young's inequality give
\[
 \left\|\left(\sum_\theta|\chi_j(D_x)F_\theta|^2\right)^{1/2}\right\|_p
 \le \|\check\chi_j\|_1
 \left\|\left(\sum_\theta|F_\theta|^2\right)^{1/2}\right\|_p.
\]
Only finitely many such terms occur. A fixed dilation replaces $R$ by a
fixed multiple of $R$, which is harmless in the displayed logarithms.

\begin{proof}[Proof of Theorem~\ref{cor:elliptic}]
We verify the frequency inclusions for three model surfaces.

\emph{The sphere.} Let $\Sigma=\{\eta:|\eta|=\lambda\}$, where $\lambda>0$
is fixed. Its caps may be described by
\[
 \bigl||\eta|-\lambda\bigr|\le C\delta,
 \qquad \left|\frac\eta{|\eta|}-\omega_\theta\right|\le C\sqrt\delta.
\]
For $|s-\lambda|\le c\delta$,
\[
 |s-|\eta||\le |s-\lambda|+\bigl||\eta|-\lambda\bigr|
 \le(C+c)\delta.
\]
Thus the product cap is contained in a fixed enlargement of the cone
sector about $\omega_\theta$. Its radial coordinate remains in a compact
positive interval. The tensor lemma applies. A frequency translation and
an invertible linear transformation give every ellipsoid.

\emph{The paraboloid.} Consider
\[
 \Sigma=\{(u,v):v=|u|^2,\ |u|\le M\},\qquad u\in\R^{n-1},
\]
where $M$ is fixed. Introduce one extra frequency $s$, and make the
invertible linear change
\begin{equation}\label{eq:paraboloidlift}
 \xi'=u,\qquad \xi_n=\frac{v-s}{2},\qquad
 \tau=\frac{v+s}{2}.
\end{equation}
The slice $s=1$ of the positive cone is precisely $v=|u|^2$. Indeed,
\[
 \tau^2-|\xi|^2=sv-|u|^2.
\]
If $|v-|u|^2|\le C\delta$ and $|s-1|\le c\delta$, then
\[
 |\tau^2-|\xi|^2|
 \le |v-|u|^2|+|s-1|\,|v|\lesssim_{M,C,c}\delta.
\]
For sufficiently small $\delta$, $\tau$ and $|\xi|$ remain bounded above
and bounded away from zero. Division by $\tau+|\xi|$ therefore gives
$|\tau-|\xi||\lesssim\delta$.

The angular map of the central slice is
\[
 \Omega(u)=\frac{(2u,|u|^2-1)}{1+|u|^2}.
\]
It is uniformly bi-Lipschitz on the fixed parameter ball, since
\[
 |\Omega(u)-\Omega(w)|
 =\frac{2|u-w|}{\sqrt{(1+|u|^2)(1+|w|^2)}}.
\]
On the lifted product cap,
\[
 \left|\frac\xi{|\xi|}-\Omega(u_\theta)\right|
 \lesssim_M |u-u_\theta|+|v-|u|^2|+|s-1|
 \lesssim\sqrt\delta.
\]
Consequently every product cap lies in an enlarged canonical cone sector.
The separated centers $u_\theta$ have separated images at the same scale,
up to fixed constants. Matching them to the nearest centers of a
canonical angular net has bounded multiplicity. Apply the tensor lemma
with the transformation \eqref{eq:paraboloidlift} and the translation
corresponding to $s=1$.

Every graph
\[
 v=c+\ell\cdot u+\tfrac12 u^T H u,
\]
with $H$ positive definite, is a fixed invertible affine image of this
model; negative definite $H$ is handled by reflection. This proves the
assertion for elliptic paraboloids.

\emph{The hyperboloid.} Consider the upper sheet
\[
 \Sigma=\{(u,v):v=\sqrt{1+|u|^2},\ |u|\le M\}.
\]
Use the extra frequency $s$ as the last spatial cone coordinate:
\[
 (\xi',\xi_n,\tau)=(u,s,v),\qquad |s-1|\le c\delta.
\]
On a product cap,
\begin{align*}
 \left|\tau-\sqrt{|\xi'|^2+\xi_n^2}\right|
 &\le\left|v-\sqrt{1+|u|^2}\right|
       +\left|\sqrt{1+|u|^2}-\sqrt{s^2+|u|^2}\right|\\
 &\lesssim_{C,c}\delta.
\end{align*}
The second term is bounded by
$|1-s^2|/(\sqrt{1+|u|^2}+\sqrt{s^2+|u|^2})$.
The cone radius lies in a fixed compact positive interval.

Here the central angular map is
\[
 \Omega(u)=\frac{(u,1)}{\sqrt{1+|u|^2}}.
\]
It and its inverse $u=\Omega'/\Omega_n$ are Lipschitz on the fixed compact
patch, because $\Omega_n\ge(1+M^2)^{-1/2}$. The product cap therefore has
angular aperture $O(\sqrt\delta)$, and matching to cone sectors again has
bounded multiplicity. The tensor lemma applies. Reflection in $v$ gives
the lower sheet, and fixed invertible affine changes give the stated
elliptic hyperboloids.

For completeness, an affine frequency change $\eta=A\zeta+a$ corresponds
to the physical transformation $g(x)\mapsto e^{ia\cdot x}g(A^Tx)$.
Both norms acquire the same factor $|\det A|^{-1/p}$. A fixed affine map
sends canonical caps to fixed enlargements of canonical caps and preserves
separation up to constants. Hence none of the affine changes above
introduces a scale-dependent factor. The bounded range of $R$ omitted in
the geometric checks has only a bounded number of caps and follows from
Cauchy--Schwarz.

Apply the endpoint cone estimate to obtain
\eqref{eq:ellipticintro}. Apply Corollary~\ref{cor:range} before taking the tensor slice to obtain
the stated lower-exponent logarithmic estimate. In
particular, no interpolation between frequency-supported subspaces is
needed here.
\end{proof}

If $f$ has Fourier support in the thickened surface, choose a smooth
bounded-overlap partition $\{\psi_\theta\}$ subordinate to the caps and
put $\widehat f_\theta=\psi_\theta\widehat f$. The preceding theorem then
gives the usual square-function estimate for $f$. This introduces no
multiplier error: the theorem has already been proved for arbitrary
functions supported in the enlarged caps, and
$\sum_\theta f_\theta=f$ exactly.

\begin{remark}[Scope of the elliptic deduction]\label{rem:ellipticscope}
The word elliptic here refers to the definite second fundamental form of
the stated quadrics. The argument does not establish
\eqref{eq:ellipticintro} or the lower-exponent estimate for every smooth
hypersurface with definite second fundamental form. An affine hyperplane
section of a quadratic cone is again a quadric. A deduction for arbitrary
smooth elliptic hypersurfaces would require an additional perturbation
argument that preserves the pointwise square-function norm; scalar
estimates on separate rescaled patches alone do not provide this.
\end{remark}

\section{Sharpness of the cone exponent}
\label{sec:sharpness}
Write $m=n-1=d-2$ and let $C_p(R)$ be the smallest constant in the reverse square-function estimate for the canonical sectors at scale $R$. The following lower bounds are independent of the external local-smoothing input.

\begin{proposition}[Necessary exponent and endpoint loss]\label{prop:sharpness}
At the proposed endpoint,
\begin{equation}\label{eq:lognecessary}
 C_{\pn}(R)\gtrsim_n[\log R]^{1/\pn},
 \qquad \pn=\frac{2(m+1)}m.
\end{equation}
For every $p>\pn$, there is $\sigma_{n,p}>0$ such that
$C_p(R)\gtrsim_{n,p}R^{\sigma_{n,p}}$. Thus no estimate with an arbitrarily small power loss can hold above the critical exponent. The logarithmic endpoint lower bound is compatible with the upper power $5/p_n$ in Theorem~\ref{thm:main}; it does not show that this upper power is optimal.
\end{proposition}

\subsection{The positive bush estimate}
Let $\{\omega\}\subset S^m$ be an $R^{-1/2}$-separated net of cardinality $N\asymp R^{m/2}$ and define dual planks
\[
 P_\omega=\{(x,t):|t+x\cdot\omega|\le C,\quad
 |x-(x\cdot\omega)\omega|\le C\sqrt R,\quad
 |x\cdot\omega-t|\le CR\},\qquad M_R=\sum_\omega\1_{P_\omega}.
\]
Their volume is comparable to $R^{(m+2)/2}$. We claim
\begin{equation}\label{eq:bushmoment}
 \int M_R^s\lesssim_{m,s}R^{m+1}
\end{equation}
for $1\le s<m/(m-1)$ if $m\ge2$, and for $s=2$ or $2<s<4$ if $m=1$. In particular the claim holds at $s=\pn/2=(m+1)/m$ and immediately above it.

Here are the counting details. Put $r=|x|$. The planks force $r\lesssim R$ and $|t|\le r+C$. If $1\le r\le\sqrt R$ and $m\ge2$, spherical slicing gives $M_R(x,t)\lesssim N/r$. If $\sqrt R\le r\lesssim R$, put $u=r-|t|$. Then $-C\le u\lesssim R/r$ and
\[
 M_R(x,t)\lesssim N r^{-m/2}(1+\max(u,0))^{(m-2)/2}.
\]
These spherical-area bounds also apply to the separated net: enlarging a center by an angular ball of radius $cR^{-1/2}$ changes the thin radial condition by $O(1)$. For $r\le\sqrt R$ this is immediate. For $r\ge\sqrt R$, the transverse condition puts the angle within $O(\sqrt R/r)$ of one pole, so the change is
$O(r(\sqrt R/r)R^{-1/2}+r/R)=O(1)$. Disjoint angular balls therefore justify counting by spherical area even in the narrow strips.

For $m\ge2$, the inner contribution is at most
\[
 N^s\int_1^{\sqrt R}r^{m+1-s}\dd r
 \lesssim R^{(m-1)s/2+(m+2)/2}\le R^{m+1}.
\]
The outer contribution is at most
\[
 N^sR^{1+(m-2)s/2}
 \int_{\sqrt R}^{CR}r^{m-1-(m-1)s}\dd r
 \lesssim R^{m+1},
\]
since $s<m/(m-1)$. The bounded region $r\le1$ is smaller. For $m=1$, spherical slicing gives instead
$M_R\lesssim Nr^{-1/2}(1+\max(u,0))^{-1/2}$. If $2<s<4$, the $u$ integral is bounded and the outer contribution is
$N^s\int_{\sqrt R}^{CR}r^{1-s/2}\dd r\lesssim R^2$; the inner one is
$O(N^sR^{1-s/4})\le CR^2$. At $s=2$, the $u$ integral introduces logarithms. The outer contribution is
$N^2\int_{\sqrt R}^{CR}(1+\log(CR/r))\dd r\lesssim R^2$, and the inner one is $O(R^{3/2}\log R)\lesssim R^2$. This proves \eqref{eq:bushmoment}.

The same estimate holds for uniformly Schwartz canonical packet envelopes: bound each squared packet by
$C_L\sum_{k\ge0}2^{-Lk}\1_{2^kP_\omega}$ and use Minkowski and isotropic dilation, choosing $L$ larger than the dimension-dependent dilation exponent.

\subsection{Polynomial failure above the endpoint}
Take one smooth Fourier-supported packet in each sector, normalized to size one near the origin and with a nonnegative Fourier profile. On a sufficiently small fixed ball their phases lie in a common half-plane, so
$\|\sum_\omega f_\omega\|_p\gtrsim N$. For $p>\pn$ sufficiently close to $\pn$, \eqref{eq:bushmoment} gives
\[
 \left\|\left(\sum_\omega|f_\omega|^2\right)^{1/2}\right\|_p
 \lesssim R^{(m+1)/p},\qquad
 C_p(R)\gtrsim R^{m/2-(m+1)/p}.
\]
This is a positive power. To cover every larger exponent, fix one such $r_*>\pn$ with $r_*\le p$. The same square function $S$ has $\|S\|_\infty\lesssim\sqrt N$, so ordinary scalar interpolation gives
\[
 \|S\|_p\lesssim
 R^{m(1-r_*/p)/4+(m+1)/p},\qquad
 C_p(R)\gtrsim R^{m/4+(mr_*/4-m-1)/p}.
\]
The last exponent is positive for all $p\ge r_*>2(m+1)/m$. No interpolation on Fourier-supported subspaces is invoked.

\subsection{A smooth endpoint example}
Choose nonnegative nonzero $a\in C_c^\infty((1,2))$ and even nonnegative nonzero $\chi\in C_c^\infty((-c,c))$. With $\delta=R^{-1}$, set
\[
 \widehat f_\delta(\xi,\tau)=\delta^{-(m+2)/2}a(|\xi|)
                 \chi\bigl((\tau-|\xi|)/\delta\bigr).
\]
Split by a smooth angular partition at aperture $\sqrt\delta$. Each piece has a uniformly bounded rescaled symbol and hence a unit-amplitude Schwartz envelope of the preceding type. Thus
\begin{equation}\label{eq:endpointsquareupper}
 \left\|\left(\sum_\theta|f_\theta|^2\right)^{1/2}\right\|_{\pn}^{\pn}
 \lesssim R^{m+1}.
\end{equation}
Integration in $\tau$ gives
\[
 f_\delta(x,t)=C R^{m/2}\check\chi(t/R)W(x,t),\qquad
 W(x,t)=\int_{\R^{m+1}}e^{i(x\cdot\xi+t|\xi|)}a(|\xi|)\dd\xi.
\]
Stationary phase on $S^m$, uniformly for large $r=|x|$, gives
\[
 W(x,t)=r^{-m/2}\{c_+A(t+r)+c_-A(t-r)\}+O(r^{-m/2-1}),
 \quad A(s)=\int a(\rho)\rho^{m/2}e^{i\rho s}\dd\rho,
\]
with nonzero constants $c_\pm$. Since $A(0)>0$ and $A$ decays rapidly, the first term dominates when $|t+r|\le c_1$ and $r\ge r_0$, for suitable fixed constants. Taking $c_0>0$ small also keeps $|\check\chi(t/R)|$ bounded below for $r\le c_0R$. Therefore
\[
 |f_\delta(x,t)|\gtrsim R^{m/2}r^{-m/2}
 \quad(r_0\le r\le c_0R,\ |t+r|\le c_1).
\]
The shell volume element is comparable to $r^m\dd r\dd t$, so
\[
 \|f_\delta\|_{\pn}^{\pn}
 \gtrsim R^{m\pn/2}\int_{r_0}^{c_0R}r^{m-m\pn/2}\dd r
 \asymp R^{m+1}\log R.
\]
Together with \eqref{eq:endpointsquareupper}, this proves \eqref{eq:lognecessary} and completes Proposition~\ref{prop:sharpness}.

\section{Acknowledgements}

AI tools were used in all aspects of producing this manuscript, including discovery, writing, checking, formalization etc.

\end{document}